\documentclass{gtart}

\def\ifplaintex{\expandafter\ifx\csname documentclass\endcsname\relax}

\ifplaintex 
\else
\fi

\expandafter\ifx\csname beginpicture\endcsname\relax
\expandafter\ifx\csname documentclass\endcsname\relax
\input pictex \else
\input prepictex \input pictex \input postpictex \fi\fi

\def\gt{{\mathsurround=0pt\it $\cal G\mskip-2mu$eometry \&\ 
$\cal T\!\!$opology}}        %  journal title in recommended style

\def\gtp{{\mathsurround=0pt\it $\cal G\mskip-2mu$eometry \&\ 
$\cal T\!\!$opology $\cal P\!$ublications}}  % GT publications

\def\lognumber#1{\def\thelognumber{#1}}
\def\volumenumber#1{\def\thevolumenumber{#1}}
\def\papernumber#1{\def\thepapernumber{#1}}
\def\volumeyear#1{\def\thevolumeyear{#1}}

\def\pagenumbers#1#2{\def\startpage{#1}\def\finishpage{#2}}
\def\published#1{\def\publishdate{#1}}
\def\proposed#1{\def\theproposer{#1}}
\def\seconded#1{\def\theseconders{#1}}
\def\received#1{\def\receiveddate{#1}}

\def\accepted#1{\def\accepteddate{#1}}

\long\def\asciiabstract#1{\long\def\theasciiabstract{#1}}

\let\\\par\let\thelognumber\relax
\let\thevolumenumber\relax\let\thepapernumber\relax
\let\thevolumeyear\relax\let\thesamplenumber\relax\let\startpage\relax
\let\finishpage\relax\let\publishdate\relax\let\receiveddate\relax
\let\reviseddate\relax\let\accepteddate\relax\let\theasciititle\relax
\let\theasciiauthors\relax
\let\theasciiabstract\relax
\let\theasciiemail\relax\let\theshortauthors\relax\let\theshorttitle\relax

\long\def\maketitlep{   % start of definition of \maketitlep

\count0=\startpage

\gt\hfill      %   Journal title (top left) 
\beginpicture
\setcoordinatesystem units <0.33truein, 0.33truein> point at 2.2 0.9
\setplotsymbol ({$\cal G$})
\plotsymbolspacing=9truept
\circulararc 315 degrees from 0 1 center at 0 0
\setplotsymbol ({$\cal T$})
\circulararc 315 degrees from 1 -1 center at 1 0
\endpicture
\break
{\small\ifx\thesamplenumber\relax % sample?  
Volume \else Sample
\fi\thevolumenumber\ (\thevolumeyear)
\startpage--\finishpage\nl
Published: \publishdate}
\vglue 0.5truein plus 0.4fil minus 0.1truein

{\parskip=0pt\leftskip 0pt plus 1fil\def\\{\par\smallskip}{\ifplaintex\large
\else\Large\fi\bf\thetitle}\par\medskip}   

\vglue 0pt plus 0.1fil 

{\parskip=0pt\leftskip 0pt plus 1fil\def\\{\par}{\sc\theauthors}
\par\medskip}

\vglue 0pt plus 0.1fil 

{\small\parskip=0pt\let\newline\\
{\leftskip 0pt plus 1fil\def\\{\par}{\sl\theaddress}\par}
\expandafter\ifx\theemail\relax    % email address?
\relax\else\vglue 5pt plus 0.02fil minus 2pt\def\\{\stdspace{\rm 
and}\stdspace} 
\cl{Email:\stdspace\tt\theemail}\fi
\ifx\theurl\relax                  % URL given?
\relax\else\vglue 5pt plus 0.02fil minus 2pt\def\\{\stdspace{\rm 
and}\stdspace}
\cl{URL:\stdspace\tt\theurl}\fi\par}

\vglue 7pt plus 0.3fil minus 3pt

{\bf Abstract}
\vglue 5pt plus 0.1fil minus 2pt

\theabstract

\vglue 7pt plus 0.3fil minus 3pt

{\bf AMS Classification numbers}\quad Primary:\quad \theprimaryclass

Secondary:\quad \thesecondaryclass

\vglue 5pt plus 0.3fil minus 2pt

{\bf Keywords}\quad \thekeywords

\vglue 10pt plus 0.5fil minus 5pt

{\small  Proposed: \theproposer\hfill Received: \receiveddate\nl
Seconded: \theseconders\hfill 
\ifx\reviseddate\relax                         % paper revised?
Accepted: \accepteddate                        % no
\else
Revised: \reviseddate                          % yes
\fi}
\eject
}       %  end of definition of \maketitlep

\let\maketitlepage\maketitlep
\let\maketitle\maketitlepage

\font\phead=cmsl9 scaled 950
\font=cmsl9 scaled 1050
\font\pnum=cmbx10 scaled 913
\font\lnum=cmbx10 
\font\pfoot=cmsl9 scaled 950
\font=cmsl9 scaled 1050
\ifplaintex
\headline{\vbox to 0pt{\vskip -4.5mm\line{\small\phead\ifnum
\count0=\startpage ISSN 1364-0380 (on line)
1465-3060 (printed) \hfill {\pnum\folio}\else\ifodd\count0\def\\{ }% 
\ifx\theshorttitle\relax\thetitle\else\theshorttitle\fi\hfill{\pnum\folio}
\else\def\\{ and }{\pnum\folio}\hfill\ifx\theshortauthors\relax\theauthors
\else\theshortauthors\fi\fi\fi}\vss}}
\footline{\vbox to 0pt{\vglue 0mm\line{\small\pfoot\ifnum\count0=\startpage
\copyright\ \gtp\hfill\else
\gt, Volume \thevolumenumber\ (\thevolumeyear)\hfill\fi}\vss
}}
\else
\makeatletter
\def\@oddhead{{\small\lhead\ifnum\count0=\startpage ISSN 1364-0380 (on line)
1465-3060 (printed) \hfill {\lnum\number\count0}\else\ifodd\count0
\def\\{ }\ifx\theshorttitle\relax \thetitle \else\theshorttitle\fi\hfill
{\lnum\number\count0}\else\def\\{ and }{\lnum\number\count0}
\hfill\ifx\theshortauthors\relax 
\theauthors\else\theshortauthors\fi\fi\fi}}\def\@evenhead{\@oddhead}
\def\@oddfoot{\small\lfoot\ifnum\count0=\startpage\copyright\ \gtp\hfill\else
\gt, Volume \thevolumenumber\ (\thevolumeyear)\hfill\fi}
\def\@evenfoot{\@oddfoot}
\makeatother
\fi

\newwrite\gtoutfile
\long\gdef\makeheadfile{  %%% start of definition of \makeheadfile
{\def\\{, }\def\s{ }
\immediate\openout\gtoutfile head.xxx
\immediate\write\gtoutfile{To: math@arxiv.org}
\immediate\write\gtoutfile{Subject: put or rep NNNNN:pppp}
\immediate\write\gtoutfile{--text follows this line--}
\immediate\write\gtoutfile{Proxy-for: \ifx\theasciiauthors\relax
\theauthors\else\theasciiauthors\fi\s<\ifx\theasciiemail\relax\theemail\else\theasciiemail\fi>}
\immediate\write\gtoutfile{\noexpand\\}
\immediate\write\gtoutfile{Authors: \ifx\theasciiauthors\relax
\theauthors\else\theasciiauthors\fi}
{\def\\{ }\immediate\write\gtoutfile{Title: \ifx\theasciititle\relax
\thetitle\else\theasciititle\fi}}
\immediate\write\gtoutfile{Subj-class: GT or SG or MG etc}
\immediate\write\gtoutfile{MSC-class: \theprimaryclass\ifx\thesecondaryclass\relax\else, \thesecondaryclass\fi}
\immediate\write\gtoutfile{Journal-ref: Geom. Topol. \thevolumenumber
(\thevolumeyear) \startpage-\finishpage}
\immediate\write\gtoutfile{Comments: Published by Geometry and Topology at}
\immediate\write\gtoutfile{\s\s http://www.maths.warwick.ac.uk/gt/GTVol\thevolumenumber/paper\thepapernumber.abs.html}
\immediate\write\gtoutfile{\noexpand\\}
\immediate\write\gtoutfile{}
\ifx\theasciiabstract\relax
\immediate\write\gtoutfile{\theabstract}\else
\immediate\write\gtoutfile{\theasciiabstract}\fi
\immediate\write\gtoutfile{}
\immediate\write\gtoutfile{\noexpand\\}
\immediate\write\gtoutfile{}
\immediate\closeout\gtoutfile}}  %%% end of definition of \makeheadfile

\def\maketitlepage{\maketitlep\makeheadfile}
\let\maketitle\maketitlepage

\lognumber{325}
\volumenumber{7}\papernumber{18}
\volumeyear{2003}
\pagenumbers{641}{643}
\received{1 June 2003}
\accepted{21 September 2003}
\published{22 October 2003}
\proposed{Robion Kirby}
\seconded{Tomasz Mrowka, Cameron Gordon}

\newtheorem{theorem}{Theorem}

\def\zz{{\bf Z}}
\def\calc{\mathcal{C}}
\def\calg{\mathcal{G}}
\def\cala{\mathcal{A}}
\def\co{\colon\thinspace}

\begin{document}

\title {Splitting the concordance group\\of algebraically slice knots}
\author{Charles Livingston}
\address{Department of Mathematics, Indiana University\\Bloomington, IN  47405,
USA}
\email{livingst@indiana.edu}

\begin{abstract}{\small As a corollary of work of Ozsv\'ath and Szab\'o
\cite{os}, it is shown that the classical concordance group of
algebraically slice knots has an infinite cyclic summand and in
particular is not a divisible group.}
\end{abstract}

\asciiabstract{As a corollary of work of Ozsvath and Szabo
[Geom. Topol. 7 (2003) 615-639], it is shown that the classical
concordance group of algebraically slice knots has an infinite cyclic
summand and in particular is not a divisible group.}



\maketitlepage

Let  $\cala$ denote the concordance group of
algebraically slice knots, the kernel of Levine's homomorphism $\phi  \co
\calc \to \calg$, where $\calc$ is the classical knot concordance group and $\calg$ is
Levine's algebraic concordance group~\cite{le1}.  Little is known
about the algebraic structure of $\cala$: it is countable and
abelian, Casson and Gordon~\cite{cg1} proved that
$\cala$ is nontrivial, Jiang~\cite{j} showed it contains a subgroup
isomorphic to
$\zz^\infty$, and the author~\cite{li} proved that
it contains a subgroup isomorphic to
$\zz_2^\infty$.   We add the following theorem, a quick corollary of recent work of 
Ozsv\'ath and Szab\'o~\cite{os}.

 \begin{theorem} The group $\cala$ contains a summand isomorphic to
$\zz$ and in particular  $\cala$ is not divisible.\end{theorem}

\begin{proof} In~\cite{os} a homomorphism  $\tau \co \calc
\to
\zz$ is constructed.  We prove that $\tau$ is
nontrivial on
$\cala$. The theorem follows since, because Im$(\tau)$ is free,  there is
the induced splitting,
$\cala
\cong
\mbox{Im}(\tau) \oplus  \mbox{Ker}(\tau)$.  No element representing a
generator of $\mbox{Im}(\tau)$ is divisible.

According to~\cite{os},  $|\tau(K)| \le g_4(K)$, where $g_4$ is the 4--ball genus of a knot, and there is the example of the $(4,5)$--torus knot  $T$  for
which $\tau(T) = 6$.  We will show
that there is a knot  $T^*$ algebraically concordant to $T$  with
$g_4(T^*) < 6$.  Hence, $T \# - T^*$ is an
algebraically slice knot with nontrivial $\tau$, as desired.

Recall that $T$ is a fibered knot with fiber $F$ of genus $(4-1)(5-1)/2
=6$. Let   $V$ be the $12 \times 12$
Seifert matrix for $T$ with respect to  some basis for $H_1(F)$.  The
quadratic form $q(x) = x V x^t$ on
$\zz^{12}$ is equal to the form given by
$(V + V^t)/2$.  Using~\cite{glm} the signature of this symmetric bilinear form can be computed to be  8, so $q$ is indefinite, and thus
by Meyer's theorem~\cite{hm} there is a nontrivial  primitive element
$z$ with $q(z) = 0$.  Since $z$ is primitive, it is a  member   of a
symplectic basis for $H_1(F)$. Let $V^*$ be the Seifert matrix for $T$
with respect to that basis.  The canonical construction of a
Seifert surface  with Seifert matrix $V^*$ (\cite{se}, or see~\cite{bz})  yields a surface $F^*$ such that
$z$ is represented by a simple closed curve on
$F^*$ that is unknotted in $S^3$.  Hence,   $F^*$ can be surgered in the $4$--ball to show that its boundary $T^*$ satisfies 
$g_4(T^*) < 6$.  Since $T^*$ and $T$ have the same Seifert form, they are
algebraically concordant.  
\end{proof}

{\bf Addendum}\qua An alternative proof of Theorem 1 follows from 
the construction of knots with trivial Alexander polynomial for which 
$\tau$ is nontrivial, to appear in a forthcoming paper.


\begin{thebibliography}

\bibitem{bz} {\bf G Burde}, {\bf    H Zieschang},  
{\it Knots},
 de Gruyter Studies in Mathematics, 5, Walter de Gruyter \& Co., Berlin (1985) 

\bibitem{cg1} {\bf A Casson}, {\bf C McA Gordon}, {\it Cobordism of
classical knots}, from: ``A la recherche de la Topologie perdue'',
(Guillou and Marin, editors), Progress in Mathematics, Volume 62
(1986), originally published as an Orsay Preprint (1975)


\bibitem{glm} {\bf C  McA  Gordon}, {\bf  R   Litherland},  {\bf K Murasugi}, 
 {\it Signatures of covering links}, 
  Canad. J. Math.  33 (1981)   381--394 

\bibitem{hm} {\bf D  Husemoller},  {\bf J Milnor},
  {\em Symmetric Bilinear Forms}.
  Ergebnisse der Mathematik und ihrer Grenzgebiete, Band 73,
Springer-Verlag, New York-Heidelberg  (1973)

\bibitem{j} {\bf B Jiang},
 {\it  A simple proof that the concordance group of algebraically
 slice knots is infinitely generated},
  Proc. Amer. Math. Soc. 83  (1981)  189--192 



\bibitem{le1} {\bf J Levine}
{\it Knot cobordism groups in codimension two}, 
  Comment. Math. Helv. 44  (1969)  229--244  

 \bibitem{li} {\bf C  Livingston},
{\it  Order 2 algebraically slice knots}, from 
  Proceedings of the Kirbyfest   (Berkeley, CA, 1998)  335--342, Geom. Topol. Monogr., 2,
Geom. Topol. Publ., Coventry (1999)


\bibitem{os} {\bf P Ozsv\'ath}, {\bf Z Szab\'o}, {\it Knot Floer
homology and the four-ball genus}, Geometry and Topology 7 (2003)
615-639, {\tt arXiv:math.GT/0301149}


\bibitem{se} {\bf H Seifert},
{\it   \"{U}ber das Geschlecht von Knoten}, 
  Math. Ann.  110 (1934) 571--592 



\end{thebibliography}
\end{document}